\documentclass{article}
\usepackage{amssymb}

\usepackage{graphicx}
\usepackage{amsmath}

\newtheorem{theorem}{Théorème}[section]

\newtheorem{corollary}[theorem]{Corollaire}

\newtheorem{proposition}[theorem]{Proposition}

\newenvironment{proof}[1][Preuve]{\textbf{#1.} }{\ \rule{0.5em}{0.5em}}
\input{tcilatex}

\begin{document}

\title{{\Large STRUCTURE D'UN DRAPEAU\ RIEMANNIEN}}
\author{\textit{Hassimiou Diallo} \\
%EndAName
\textit{Laboratoire de Math\'{e}matiques Fondamentales, }\\
\textit{ENS, Universit\'{e} de Cocody}\\
\textit{08 B.P. 10 ,\ Abidjan, C\^{O}TE D'IVOIRE.}\\
\textit{Email:diallomh@yahoo.fr}}
\maketitle

\begin{abstract}
\textit{The aim of this paper is to show that any stable complete Riemannian
flag on a compact and connected manifold is conjugated to a flag of
homogenus foliations ( see D\'{e}finitions).}

\textit{\ Also, we give a characterization of complete Riemannian flags that
are homogenus. }

This result is a step toward the classification of Riemannian flags.
\end{abstract}

\section{Introduction\newline
}

L'objet de ce travail \ est l'\'{e}tude d'une classe de drapeaux riemanniens
. Avant d'exposer la position du probl\`{e}me, il est utile de rappeller les
d\'{e}finitions suivantes.

\textit{Quelques d\'{e}finitions:}

1) Un \textit{drapeau de feuilletages} (ou tout simplement un \textit{drapeau%
}) sur une va-ri\'{e}t\'{e} $M$ est la donn\'{e}e sur cette vari\'{e}t\'{e} $%
\ $d'une suite croissante de feuilletages $\mathcal{D=}$($\mathcal{F}%
_{1},...,\mathcal{F}_{p}$) telle que $\dim \mathcal{F}_{i}=i$. L'entier $p$
est la longueur du drapeau. Le drapeau est dit \textit{complet}\textbf{\ }si
le dernier feuilletage est de codimension \textit{un}. Un drapeau de longueur%
\textit{\ un} est r\'{e}duit \`{a} son flot.

2) Un drapeau\textbf{\ }est dit\textbf{\ }\textit{riemannien}\textbf{\ }( 
\textit{de Lie}$,$\textit{\ minimal}$,$\textit{\ parall\'{e}lisable}) si
chacun de ses feuilletages est riemannien (de Lie, minimal, transversalement
parall\'{e}lisable).

3) Un drapeau sur $M=\mathbb{R}^{n}$ ou $\mathbb{T}^{n}$ est dit \textit{%
lin\'{e}aire} si chacun de ses feuilletages est lin\'{e}aire.

4) Etant donn\'{e}s un drapeau $\mathcal{D}=(\mathcal{F}_{1},\ldots ,%
\mathcal{F}_{p}$ $)$ sur une vari\'{e}t\'{e} $M$ et une submersion $\pi $
d'une vari\'{e}t\'{e} $M^{\#}$ sur $M$, un \textit{rel\`{e}vement du drapeau}
$\mathcal{D}$ \ sur $M^{\#}$ est, s'il existe, un drapeau $\mathcal{D}^{\#}=(%
\mathcal{F}_{1}^{\#},\ldots ,\mathcal{F}_{p}^{\#}$ ) sur $M^{\#}$, tel que
pour tout $i\in \{1,\ldots ,$ $p-1\}$, le feuilletage $\mathcal{F}_{i}^{\#}$
\ se projette sur le feuilletage $\mathcal{F}_{i}$ \ par $\pi $.

5) Un drapeau $\mathcal{D}=\left\{ \mathcal{F}_{1},\ldots ,\mathcal{F}%
_{p}\right\} $ sur une vari\'{e}t\'{e} $M$ et un drapeau \newline
$\mathcal{D}^{\prime }=(\mathcal{F}_{1}^{\prime },\ldots ,\mathcal{F}%
_{p}^{\prime }$ $)$ sur une vari\'{e}t\'{e} $M^{\prime }$ sont dits \textit{%
conjugu\'{e}s} s'il existe un diff\'{e}omorphisme de $M$\ sur $M^{\prime }$
\'{e}changeant $\mathcal{F}_{i}$ et $\mathcal{F}_{i}^{\prime }$.

6) Un drapeau de feuilletages riemanniens $\mathcal{D}$ est dit \textit{%
stable} si tout feuilletage du drapeau est soit minimal, soit r\'{e}gulier
et d'adh\'{e}rence dans $\mathcal{D}$ .

Si $F$ est l'adh\'{e}rence d'une feuille d'un feuilletage $\mathcal{F}_{r}$
\ d'un drapeau stable $(M,\mathcal{D=}(\mathcal{F}_{1},...,\mathcal{F}%
_{p})), $ les restrictions des feuilletages $\mathcal{F}_{1},...,\mathcal{F}%
_{s}$ \`{a} $F$, o\`{u} $r\leq s<\min (p,\dim F)$, \ forment un drapeau
qu'on appellera\textit{\ trace\ }de $\mathcal{D}$ sur $F$ .

\textit{Position du probl\`{e}me}

L'objectif dans l'\'{e}tude des drapeaux de feuilletages, a \'{e}t\'{e} au d%
\'{e}part en partant des r\'{e}sultats de Y. Carri\`{e}re, de E. F\'{e}dida
et de P. Molino sur les feuilletages(\cite{CAR}, \cite{FED}, \cite{MOL2} ) ,
de voir ce qu'on peut dire des drapeaux de feuilletages de Lie et des
drapeaux riemanniens en g\'{e}n\'{e}ral.

Les r\'{e}sultats obtenus sont les suivants:

Sur une vari\'{e}t\'{e} compacte,

1) Tout drapeau riemannien minimal est conjugu\'{e} \`{a} un drapeau
lin\'{e}aire \cite{BD}.

2) Tout drapeau riemanien se rel\`{e}ve sur le fibr\'{e} des rep\`{e}res
orthonorm\'{e}s directs transverses du dernier feuilletage en un drapeau de
feuilletages transversalement parall\'{e}lisables. Ce qui a permis
d'\'{e}tablir qu'un drapeau riemannien complet d'une \textit{vari\'{e}t\'{e}
compacte} est un drapeau parall\'{e}lisable\cite{DIA2}.

3) Le th\'{e}or\`{e}me de rel\`{e}vement de E. F\'{e}dida \ permet
d'\'{e}crire un th\'{e}or\`{e}me de structure des drapeaux de Lie\cite{DIA1}.

Ici, il s'agit d'aborder l'\'{e}tude de la structure d'un type particulier
de drapeau riemannien.

De fa\c{c}on pr\'{e}cise ,

1) on montre que \textit{tout drapeau riemannien complet stable d'une
vari\'{e}t\'{e} compacte est conjugu\'{e} \`{a} un drapeau de feuilletages
homog\`{e}nes.}

2) On \'{e}tablit \textit{une condition n\'{e}cessaire et suffisante pour
qu'un drapeau riemannien complet d'une vari\'{e}t\'{e} compacte soit un
drapeau de feuilletages homog\`{e}nes}.

Le 1) permet alors de statuer sur les drapeaux riemanniens complets ou non
d'une \textit{3-vari\'{e}t\'{e}} \textit{compacte }et sur les drapeaux
riemanniens complets d'une\textit{\ 4-vari\'{e}t\'{e} compacte, }et
r\'{e}pr\'{e}sente un pas important vers la classification des drapeaux
riemanniens.

Il s'agira pour la suite d' achever la classification de ces drapeaux
lorsqu'ils sont support\'{e}s par une \textit{n-}vari\'{e}t\'{e} compacte ($%
n\geq 4),$ qu'ils soient complets ou non , stables ou non.

Dans tout ce qui suit les objets sont $C^{\infty }$, \ la vari\'{e}t\'{e} $M$
est connexe et orientable, et les feuilletages consid\'{e}r\'{e}s sont
suppos\'{e}s orientables.\newpage

\section{Structure d'un drapeau riemannien complet}

Pour \'{e}tablir le th\'{e}or\`{e}me de struture d'un drapeau riemannien
complet, nous allons montrer que pour un drapeau riemannien , admet toujours
une m\'{e}trique quasi-fibr\'{e}e par rapport \`{a} chaque feuilletage du
drapeau.

\begin{proposition}
\bigskip \label{prop.metrique qf}Etant donn\'{e} un drapeau riemannien , il
existe toujours une m\'{e}trique quasi-fibr\'{e}e commune \`{a} tous les
feuilletages.
\end{proposition}

\begin{proof}
\bigskip\ Il suffit de prouver que si $\mathcal{F}$ et $\mathcal{F}^{\prime }
$ sont deux feuilletages riemanniens d'une vari\'{e}t\'{e} $M$ tels que $%
\mathcal{F}^{^{\prime }}$soit une extension de $\mathcal{F}$, alors , il
existe sur $M$ une m\'{e}trique riemannienne quasi-fibr\'{e}e \`{a} la fois
pour $\mathcal{F}$ et $\mathcal{F}^{\prime }.$ Et on ach\`{e}ve par r\'{e}%
currence.

Soit $\left( U_{i}\right) _{i\in I}$ un recouvrement de $M$ \ form\'{e}
d'ouverts distingu\'{e}s \`{a} la fois pour $\mathcal{F}$ et $\ \mathcal{F}%
^{\prime }$ de sorte que les feuilletages riemanniens soient d\'{e}finis
respectivement par les cocycles $\left( U_{i},T,f_{i},g_{ij}\right) _{i\in I}
$ et $\left( U_{i},T^{\prime },\text{ }f_{i}^{\prime },g_{ij}^{\prime
}\right) _{i\in I}$ . Pour chaque $i\in I$, \ il existe une submersion
locale $\theta _{i}$ de $f_{i}$($U_{i})$ dans $T^{\prime }$ telle $%
f_{i}^{\prime }=\theta _{i}\circ f_{i}.$ Si $\ U_{i}\cap U_{j}\neq \emptyset 
$ , on a le diagramme commutatif \ suivant $:$%
\begin{equation*}
\begin{array}{ccccc}
U_{i}\cap U_{j} & \overset{f_{i}}{\rightarrow } & T & \overset{\theta _{i}}{%
\rightarrow } & T^{\prime } \\ 
^{Id_{U_{i}\cap U_{j}}}\downarrow  &  & \downarrow ^{g_{ji}} &  & \downarrow
^{g_{ji}^{\prime }} \\ 
U_{i}\cap U_{j} & \overset{f_{j}}{\rightarrow } & T & \overset{\theta _{j}}{%
\rightarrow } & T^{\prime }
\end{array}
\end{equation*}

Et si $h$ et $h^{\prime }$ sont les m\'{e}triques riemanniennes respectives
de $T$ et $T^{\prime }$d\'{e}finissant les feuilletages riemanniens $%
\mathcal{F}$ et $\ \mathcal{F}^{\prime },$ comme $g_{ji}^{\prime }\circ
\theta _{i}=\theta _{j}\circ g_{ji}$ , et que $g_{ji}$ est une isom\'{e}trie
locale, alors 
\begin{equation*}
Tg_{ji}(KerT\theta _{i})=KerT\theta _{j}\text{ \ et \ }Tg_{ji}((KerT\theta
_{i})^{\bot })=(KerT\theta _{j})^{\bot }
\end{equation*}
\newline
Sur chaque ouvert $\ f_{i}$($U_{i})$ modifions la m\'{e}trique induite par $h
$ de la fa\c{c}on suivante: tout en la pr\'{e}servant sur $KerT\theta _{i}$
et en gardant l'orthogonalit\'{e} des sous-fibr\'{e}s $KerT\theta _{i}$ et ($%
KerT\theta _{j})^{\bot },$ on la substitue sur $(KerT\theta _{i}$)$^{\bot }$
par l$^{\prime }$image r\'{e}ciproque par $\theta _{i}$ de la m\'{e}trique
riemannienne $h^{\prime }$de $T^{\prime }$ $.$\ Avec ces nouvelles m\'{e}%
triques sur les $f_{i}$($U_{i})$, les $\theta _{i}$ deviennent par
construction des submersions riemanniennes, et les $g_{ij}$ restent encore
des isom\'{e}tries . Comme la vari\'{e}t\'{e} transverse mod\`{e}le $T$ du
feuilletage riemanien $\mathcal{F},$ est par construction , la somme
disjointe des $f_{i}$($U_{i})$ ,\ alors la somme disjointe de ces nouvelles m%
\'{e}triques sur les $f_{i}$($U_{i})$ d\'{e}finissent sur $T$ une m\'{e}%
trique $\widetilde{h}$ . A cette nouvelle m\'{e}trique de $T,$ il est
loisible d'y associer des m\'{e}triques $\mathcal{F-}$quasi-fibr\'{e}es. Et
chacune de ces m\'{e}triques est aussi quasi-fibr\'{e}e pour $\mathcal{F}%
^{\prime }.$
\end{proof}

Ainsi , \'{e}tant donn\'{e} un drapeau riemannien, il existe toujours une
m\'{e}trique quasi-fibr\'{e}e commune \`{a} tous les feuilletages. On dira
d'une telle m\'{e}trique qu'elle est \textit{quasi-fibr\'{e}e pour le drapeau%
}.

Ceci \'{e}tant , le \ premier r\'{e}sultat porte sur la structure d'un
drapeau riemannien complet stable.

Nous rappelons qu'un drapeau de feuilletages riemanniens $\mathcal{D}$ est
dit \textit{stable} si tout feuilletage du drapeau est soit \textit{minimal}%
, soit \textit{r\'{e}gulier et d'adh\'{e}rence dans} $\mathcal{D}$ .Cela
signifie que chaque feuilletage \`{a} feuilles non ferm\'{e}es est dense
dans les feuilletages suivants (\cite{MOL2}) jusqu'au prochain feuilletage
\`{a} feuilles ferm\'{e}es \ du drapeau s'il en contient.

\subsection{Structure d'un drapeau riemannien complet stable}

\begin{proposition}
\label{proposition2} Soit $(M,\mathcal{F}_{1,}..,\mathcal{F}_{k+1})$ un
drapeau d'une $n+1-$vari\'{e}t\'{e} compacte connexe tel que:

i) le flot est un feuilletage transversalement parall\'{e}lisable,

ii) la vari\'{e}t\'{e} basique de ce flot admet pour rev\^{e}tement
universel un groupe de Lie\ r\'{e}soluble,

iii)\ le feuilletage $\mathcal{F}_{k+1}$ co\"{i}ncide avec le feuilletage
basique du flot.

Alors la vari\'{e}t\'{e} $M$ est diff\'{e}omorphe \`{a} une vari\'{e}t\'{e}
homog\`{e}ne d'un groupe de Lie r\'{e}soluble et le drapeau est conjugu\'{e} 
\`{a} un drapeau de feuilletages homog\`{e}nes.
\end{proposition}

\begin{proof}
\bigskip L'entier $k$ \'{e}tant la dimension de l$^{\prime }$a$\lg $\`{e}bre
de Lie structurale du flot $(M,\mathcal{F)}$ , si $\mathbb{T}%
^{k+1}\hookrightarrow M\overset{\pi _{k}}{\rightarrow }W$ est sa fibration
basique , alors, on sait qu'il existe une matrice $A\in \mathbb{SL(}k+1;%
\mathbb{Z)}$ diagonalisable ayant toutes ses valeurs propres positives telle
que la restriction du flot \`{a} chaque fibre soit conjugu\'{e}e \`{a} un m%
\^{e}me flot mod\`{e}le lin\'{e}aire du tore $\mathbb{T}^{k+1}$de direction $%
v_{1}$ un des vecteurs propres de $A$\cite{AM}. Des exemples de telles
matrices peuvent \^{e}tre trouv\'{e}s dans\cite{EN}$.$

Puisque chacun des feuilletages contient le flot, alors il laisse invariant
les feuilletages et ensuite comme l'action de ce flot est transitive, le
drapeau induit sur chacune des adh\'{e}rences du flot$-$ le tore $\mathbb{T}%
^{k+1}-,$ le m\^{e}me drapeau $\mathcal{D}_{k}^{\prime }\mathcal{=}$($%
\mathcal{F}_{1}^{\prime },...,\mathcal{F}_{k}^{\prime }$). Comme les
feuilletages de $\mathcal{D}_{k}^{\prime }$ en tant que restriction de
feuilletages riemanniens sur une partie satur\'{e}e pour chacun des
feuilletages $\mathcal{F}_{1,}..,\mathcal{F}_{k}$ \ , sont aussi des
feuilletages riemanniens, alors d'apr\`{e}s \cite{BD}le drapeau $\mathcal{D}%
_{k}^{\prime }$ est lin\'{e}aire.

On peut pr\'{e}ciser que, si $v_{1},.....,v_{k+1}$sont les $(k+1)-$ vecteurs
propres associ\'{e}s aux $(k+1)-$ valeurs propres non n\'{e}cessairement
distinctes de $A,$le drapeau mod\`{e}le $\mathcal{D}_{k}^{\prime }$ induit
sur chaque adh\'{e}rence est tel que pour tout $i,$ 1$\leq i\leq k,$ $%
\overline{<v_{1},\ldots ,v_{i}>}$ soit une direction du feuilletage lin\'{e}%
aire $\mathcal{F}_{i}^{\prime }.$

On peut noter pour la suite que la projection lin\'{e}aire $\sigma _{i}$ de $%
\mathbb{R}^{k+1}$ sur le sous-espace $\overline{<v_{i+1},\ldots ,v_{k+1}>}%
\cong \mathbb{R}^{k-i+1}$ suivant la direction $\overline{<v_{1},\ldots
,v_{i}>}$ est une d\'{e}v\'{e}loppante de $\mathcal{F}_{i}^{\prime }$ et
qu'ensuite, puisque $A$ laisse invariant le feuilletage relev\'{e} de $%
\mathcal{F}_{i}^{\prime }$ sur $\mathbb{R}^{k+1}$et le sous-espace $%
\overline{<v_{i+1},\ldots ,v_{k+1}>}$, il existe une matrice carr\'{e}e $%
A_{i}$ d'ordre $k-i+1$, dont l'application lin\'{e}aire $\rho _{A_{i}}$
associ\'{e}e est telle qu'on ait $\rho _{A_{i}}\circ \sigma _{i}=\sigma
_{i}\circ \rho _{A}$. Il en r\'{e}sulte que cette matrice $A_{i}$ est
diagonalisable , \`{a} valeurs propres positives et de vecteurs propres $%
v_{i+1},\ldots ,v_{k}$, de sorte que pour $t\in \mathbb{R}$, la matrice $%
A_{i}^{t}$ est d\'{e}finie et l'on a: 
\begin{equation}
\rho _{A_{i}^{t}}\circ \sigma _{i}=\sigma _{i}\circ \rho _{A^{t}} 
\tag{$\star $}
\end{equation}

Ceci dit, soit $G_{W}$ \ le groupe de Lie connexe r\'{e}soluble (donc
simplement connexe \cite{DIE}) rev\^{e}tement universel de $W$ et $G=\mathbb{%
R}^{k+1}\times _{A}G_{W}$ \ la structure de groupe de Lie sur $\mathbb{R}%
^{k+1}\times G_{W}$ d\'{e}finie par la loi : 
\begin{equation*}
(x,g)\cdot (x^{\prime },g^{\prime })=(A^{\rho (g)^{-1}}x^{\prime
}+x,gg^{\prime })
\end{equation*}
o\`{u} $\rho $ est le prolongement lin\'{e}aire d'une r\'{e}pr\'{e}sention
non nulle de $\pi _{1}(W)$ sur $\mathbb{Z}$\cite{RAG}.

Ensuite, puisque la vari\'{e}t\'{e} $W=\frac{G_{W}}{\pi _{1}(W)}$ est
compacte, et $\pi _{1}(W)$ s'identifiant \`{a} un sous-groupe discret
uniforme de $G_{W}$; alors , il est clair que $\mathbb{Z}^{k+1}\times
_{A}\pi _{1}(W)$ est aussi un sous-groupe discret uniforme de $G$ et par
suite $P=\frac{G}{\mathbb{Z}^{k+1}\times _{A}\pi _{1}(W)}$ est une $n+1-$
vari\'{e}t\'{e} compacte.

Par ailleurs, comme le groupe produit semi-direct de deux groupes r\'{e}%
so-lubles est r\'{e}soluble, alors $G$ est \ r\'{e}soluble. Ensuite, le
groupe r\'{e}soluble $G_{W}$ , \'{e}tant topologiquement isomorphe \`{a} un $%
\mathbb{R}^{m}$\cite{DIE}$,$est alors contractile. Il en r\'{e}sulte que le
fibr\'{e} r\'{e}ciproque du rev\^{e}tement universel de $W$ par la
projection basique $\pi _{k}$ est triviale.

Ceci \'{e}tant, consid\'{e}rons la submersion $\Phi $ de $G=\mathbb{R}%
^{k+1}\times _{A}G_{W}$ sur $\pi _{k}^{\ast }(G_{W})$ d\'{e}finie par $\Phi
(t,g)=\psi ^{-1}(\theta (t),g$ ) o\`{u} $\theta $ est la projection
canonique de $\mathbb{R}^{k+1}$sur $\mathbb{T}^{k+1}$et $\psi :\pi
_{k}^{\ast }(G_{W})\rightarrow \mathbb{T}^{k+1}\times G_{W}$ \ une
trivialisation de $\pi _{k}^{\ast }(G_{W})$. Pour tous $s\in \pi _{1}(W)$ et
\ $(t,g)\in $ $\mathbb{R}^{k+1}\times _{A}G_{W}$ \ on a 
\begin{equation*}
\Phi (s\bullet (t,g))=s\cdot \Phi (t,g)\ 
\end{equation*}
o\`{u} la seconde action est l'action canonique de $\pi _{1}(W)$ \ sur $\pi
_{k}^{\ast }(G_{W});$ $\Phi $ passe au quotient et d\'{e}finit une
submersion surjective $\varphi $ de $P=G\diagup \mathbb{Z}^{k+1}\times
_{A}\pi _{1}(W)$ sur $M$ rendant commutatif le diagramme 
\begin{equation*}
\begin{array}{ccc}
\text{ \ \ }G & \overset{\Phi }{\rightarrow } & \pi _{k}^{\ast }(G_{W}) \\ 
^{p}\downarrow  &  & \downarrow ^{pr_{1}} \\ 
\text{ \ \ }P & \overset{\varphi }{\rightarrow } & M
\end{array}
\end{equation*}
o\`{u} $\ pr_{1}$ est la premi\`{e}re projection de $\pi _{k}^{\ast
}(G_{W})=M\times _{W}G_{W}$ \ et $p$ la projection canonique de $G$ sur $P.$
En s'aidant du diagramme commutatif suivant , 
\begin{equation*}
\begin{array}{ccccc}
\text{\ \ \ }G & \overset{\theta \times Id_{G_{W}}}{\rightarrow } & \mathbb{T%
}^{k+1}\times G_{W} & \overset{pr_{2}}{\rightarrow } & G_{W} \\ 
_{Id_{G}}\downarrow  &  & \text{ \ \ \ \ \ }\downarrow _{\psi ^{-1}} &  & 
\text{ \ \ }\downarrow _{Id_{G_{w}}} \\ 
\text{\ \ \ \ }G & \overset{\Phi }{\rightarrow } & \text{ \ \ \ \ \ \ \ }\pi
_{k}^{\ast }(G_{W}) & \overset{pr_{2}}{\rightarrow } & \text{ \ }G_{W} \\ 
_{p}\downarrow  &  & \text{ \ \ \ \ }\downarrow ^{pr_{1}} &  & \text{ \ \ }%
\downarrow _{p_{W}} \\ 
\text{\ \ \ }P & \overset{\varphi }{\rightarrow } & \text{ }M & \overset{\pi
_{k}}{\rightarrow } & \text{ }W
\end{array}
\end{equation*}
on a , pour tout $m\in M$ , 
\begin{equation*}
((\theta \times Id_{G_{W}})^{-1}\circ \psi )pr_{1}^{-1}(m)=
\end{equation*}
\begin{equation*}
(\theta \times Id_{G_{W}})^{-1}\{(z_{0},sg_{0})\in \mathbb{T}^{k+1}\times
G_{W}\diagup p_{W}(g_{0})=\pi _{k}(m)\text{ \ }et\text{ }s\in \pi _{1}(W)\}=
\end{equation*}
\begin{equation*}
(\theta \times Id_{G_{W}})^{-1}(\pi _{1}(W).(z_{0},g_{0}))=(\mathbb{Z}%
^{k+1}\times _{A}\pi _{1}(W)).(t_{0},g_{0})=p^{-1}(\varphi ^{-1}(m))\text{ ,}
\end{equation*}
o\`{u} 
\begin{equation*}
\theta (t_{0})=z_{0}=(pr_{1}\circ \psi )(m,g_{0});
\end{equation*}
alors 
\begin{equation*}
\varphi ^{-1}(m)=p((\mathbb{Z}^{k+1}\times _{A}\pi _{1}(W)).(t_{0},g_{0}))
\end{equation*}
est un singleton et il en r\'{e}sulte que $\varphi $ est une submersion
bijective, $i.e$ un isomorphisme de vari\'{e}t\'{e}s.

Le relev\'{e} horizontal du drapeau $\mathcal{D}_{k}^{\prime }$ sur $G=%
\mathbb{R}^{k+1}\times _{A}G_{W}$ \ est le relev\'{e} par $p$ du drapeau $%
\varphi ^{\ast }\mathcal{D}_{k}=$($\varphi ^{\ast }\mathcal{F}%
_{1},...,\varphi ^{\ast }\mathcal{F}_{k}).$ Ensuite en s'aidant de la
formule (*) , il est facile de voir que les feuilletages $\varphi ^{\ast }%
\mathcal{F}_{i}$ , $1\leq i\leq k,$ sont des feuilletages de Lie dont les d%
\'{e}v\'{e}loppantes respectives sont les morphismes de groupes $D_{i}$ : $%
\mathbb{R}^{k+1}\times _{A}G_{W}\rightarrow \mathbb{R}^{k}\times
_{A_{i}}G_{W}$ \ d\'{e}finis par 
\begin{equation*}
D_{i}(t,g)=(\sigma _{i}(t),g)\text{ ,}
\end{equation*}
et que le feuilletage $\varphi ^{\ast }\mathcal{F}_{k+1}$ correspond au
feuilletage d\'{e}fini par la projection $\ \mathbb{R}^{k+1}\times
_{A}G_{W}\rightarrow G_{W};$ et qu'au total le drapeau consid\'{e}r\'{e} $(M,%
\mathcal{F}_{1,}..,\mathcal{F}_{k+1})$ est \ conjugu\'{e}, au moyen de $%
\varphi ,$ \`{a} un drapeau de feuilletages homog\`{e}nes (\cite{FED}, \cite
{GHY}).
\end{proof}

\begin{theorem}
\label{theo principal} \ Tout drapeau riemannien complet stable$\mathcal{(}M,%
\mathcal{D)}$ d'une n+1-vari\'{e}t\'{e} compacte connexe est conjugu\'{e} 
\`{a} un drapeau de feuilletages homog\`{e}nes.

De fa\c{c}on plus pr\'{e}cise:

- La vari\'{e}t\'{e} $M$ est diff\'{e}omorphe \`{a} un espace homog\`{e}ne
quotient d'un groupe de Lie r\'{e}soluble $G$ par un sous-groupe discret $%
\Gamma ,$ et

- il existe une suite strictement d\'{e}croissante de sous-groupes distingu%
\'{e}s dont les orbites d\'{e}terminent le drapeau.
\end{theorem}

\begin{proof}
On pose $\dim M=n+1$ et on va raisonner par r\'{e}currence sur le nombre $l$
de feuilletages \`{a} feuilles ferm\'{e}es du drapeau.

- Si $l=0,$ le flot est minimal , le r\'{e}sultat est \'{e}vident \cite{BD}.

Sinon , soit $M_{n+1}$ une vari\'{e}t\'{e} compacte connexe supportant un
drapeau riemannien complet stable $\mathcal{D}_{n}\mathcal{=}$($\mathcal{F}%
_{1},...,\mathcal{F}_{n}$). On sait que les feuilletages du drapeau sont
tous transversalement parall\'{e}lisales\cite{DIA2}. Soit $k$ la dimension
de l'alg\`{e}bre de Lie structurale du flot de ($M_{n+1}$ ,$\mathcal{D}%
_{n});k\leq n-1$ et 
\begin{equation*}
\mathbb{T}^{k+1}\hookrightarrow M_{n+1}\overset{\pi _{k}}{\rightarrow }%
W_{n-k}
\end{equation*}
sa fibration basique.Comme le drapeau est stable et $l\neq 0$,$\overline{%
\text{ }\mathcal{F}_{1}\text{ \ }}$ , le feuilletage par les adh\'{e}rences
du flot, co\"{i}ncide avec le feuilletage $\mathcal{F}_{k+1}$ ; leurs
fibrations basiques sont alors identiques. Dans la suite excate d'Atiyah de $%
C^{0}(W_{n-k})-$ modules et d'alg\`{e}bres de Lie 
\begin{equation*}
0\rightarrow \Gamma g_{W_{n-k+1}}\rightarrow l(M,\mathcal{F}_{k+1})\overset{%
\pi _{k}}{\rightarrow }\Xi (W_{n-k})\rightarrow 0
\end{equation*}
associ\'{e}e au feuilletage \`{a} feuilles ferm\'{e}es $\mathcal{F}_{k+1}$%
\cite{MOL1}$,$ le $C^{0}(W_{n-k})-$ module\newline
$\Gamma g_{W_{n-k}}$est alors nul ; $\pi _{k}$ est donc un isomorphisme d'alg%
\`{e}bres. Le drapeau se projette sur la vari\'{e}t\'{e} basique $W_{n-k}$du
flot $\mathcal{F}_{1}$en un drapeau riemannien complet stable $\mathcal{D}%
_{n-k-1}^{\flat }\mathcal{=}$($\mathcal{F}_{1}^{\flat },...,\mathcal{F}%
_{n-k-1}^{\flat })$ o\`{u} $\ $pour tout $s,1\leq s\leq n-k-1,$ \newline
$\mathcal{F}_{s}^{\flat }$ = $\pi _{k}(\mathcal{F}_{k+1+s}).$ \ 

Comme pr\'{e}c\'{e}demment $\mathcal{D}_{n-k-1}^{\flat }$ va se projetter
sur la vari\'{e}t\'{e} basique de $\mathcal{F}_{1}^{\flat }$ en drapeau
riemannien complet stable qui est aussi le projet\'{e} de $\mathcal{D}_{n}$
sur la vari\'{e}t\'{e} basique du second feuilletage \`{a} feuilles ferm\'{e}%
es de $\mathcal{D}_{n}.$ Soit $W_{1}^{\prime }=W_{n-k}$, ..., $W_{l}^{\prime
}$, la suite des vari\'{e}t\'{e}s basiques ainsi \ obtenue par projections
successives. La projection de $\mathcal{D}_{n}$ sur $W_{l}^{\prime }$ est un
drapeau riemannien dont aucun feuilletage n'est \`{a} feuilles ferm\'{e}es
car sinon $\mathcal{D}_{n}$ contiendrait , contrairement \`{a} l'hypoth\`{e}%
se , plus de $l$ feuilletages \`{a} feuilles ferm\'{e}es. Ensuite puisque ce
drapeau sur $W_{l}^{\prime }$ est stable, ses feuilletages sont n\'{e}%
cessairement minimaux ; par suite on peut consid\'{e}rer que $W_{l}^{\prime
}=\mathbb{T}^{q}$ et que le drapeau qu'il supporte est lin\'{e}aire \cite{BD}%
. Comme le rev\^{e}tement universel $\widetilde{W_{l}^{\prime }}$ =$\mathbb{R%
}^{q}$ est un groupe de Lie r\'{e}soluble, alors par la proposition \ref
{proposition2}, la projection de $\mathcal{D}_{n}$ sur $W_{l-1}^{\prime }$
est un drapeau de feuilletages homog\`{e}nes et le rev\^{e}tement universel $%
\widetilde{W_{l-1}^{\prime }}$ est un groupe de Lie r\'{e}soluble. Par r\'{e}%
currence , toutes les vari\'{e}t\'{e}s basiques consid\'{e}r\'{e}es sont des
vari\'{e}t\'{e}s homog\`{e}nes de groupes de Lie r\'{e}solubles. Il vient
donc avec la m\^{e}me proposition \ref{proposition2}, que la vari\'{e}t\'{e} 
$M_{n+1}$ est diff\'{e}omorphe \`{a} une vari\'{e}t\'{e} homog\`{e}ne d'un
groupe r\'{e}soluble . En clair , il existe une suite d'entiers strictement
positifs $k_{1},..,k_{l-1}$ , des groupes de Lie r\'{e}solubles $%
G_{1},..,G_{l}$ , et des matrices carr\'{e}es $A_{1,..,}A_{l-1}$ tels que: 
\begin{eqnarray*}
M &\cong &\frac{G}{\pi _{1}(M)};\text{ }G=\mathbb{R}^{k_{1}}\times
_{A_{1}}G_{1}=\mathbb{R}^{k_{1}}\times _{A_{1}}(\mathbb{R}^{k_{2}}\times
_{A_{2}}G_{2})=.., \\
G_{i} &=&\mathbb{R}^{k_{i}}\times _{A_{i}}G_{i+1},G_{l-1}=\mathbb{R}%
^{k_{l-1}}\times _{A_{l-1}}\mathbb{R}^{s},G_{l}=\mathbb{R}^{q},
\end{eqnarray*}
et il est facile de voir que le rel\'{e}v\'{e} de chaque feuilletage du
drapeau $\mathcal{D}_{n}$ sur $G$ est soit le relev\'{e} d'un feuilletage
homog\`{e}ne , soit le produit du rel\'{e}v\'{e} d'un feuilletage basique et
du relev\'{e} d'un feuilletage homog\`{e}ne. Mais comme les feuilletages
basiques des feuilletages du drapeau sont aussi des feuilletages homog\`{e}%
nes, alors le drapeau $\mathcal{D}_{n}$ est \`{a} une conjugaison pr\'{e}s,
un drapeau de feuilletages homog\`{e}nes.

Enfin, la consid\'{e}ration de la suite $(H_{i})_{1\leq i\leq n}$, o\`{u} $%
H_{i}=KerD_{i}$ et les $D_{i}$ des developpantes de ces feuilletages homog%
\`{e}nes, permet de conclure.
\end{proof}

\textit{Remarques}

1) Le rev\^{e}tement universel d'une vari\'{e}t\'{e} compacte connexe
supportant un drapeau riemannien complet stable (donc un drapeau de Lie) est
un groupe de Lie finement r\'{e}soluble et le groupe d'holonomie de chacun
des feuilletages est finement polycyclique \cite{DIA1}.

2) La preuve du th\'{e}or\`{e}me\ref{theo principal} montre \ que si chacune
des matrices $A_{i}$ est r\'{e}duite \`{a} la matrice identit\'{e} alors le
groupe de Lie est ab\'{e}lien ; $G=\mathbb{R}^{n+1},$ et le drapeau est
lin\'{e}aire.

En particulier un drapeau riemannien complet dont les feuilletages sont
\`{a} feuilles ferm\'{e}es est un drapeau lin\'{e}aire du tore.

3) Les adh\'{e}rences des feuilles des feuilletages d'un drapeau stable\
sont isomorphes \`{a} des vari\'{e}t\'{e}s homog\`{e}nes de groupes de Lie
finement r\'{e}solubles. Et la trace du drapeau\ sur chacune de ces
adh\'{e}rences est un drapeau de feuilletages homog\`{e}nes.

4) La preuve de la proposition \ref{proposition2}, le th\'{e}or\`{e}me\ref
{theo principal} et \cite{DIA1} permettent de voir que tout drapeau
riemannien complet d'une $3-$vari\'{e}t\'{e} compacte connexe est , \`{a}
une conjugaison pr\`{e}s, un drapeau de Lie du tore hyperbolique $\mathbb{T}%
_{A}^{3}(trA\geq 2).$

Il en r\'{e}sulte qu'un flot riemannien sur la sph\`{e}re $\mathbb{S}^{3}$
ou sur un espace lenti-culaire de dimension $3$ n'admet aucune extension
riemannienne propre.

De m\^{e}me en dimension $\ \mathit{3}$ ou $\mathit{4},$ on a:

\begin{corollary}
Tout drapeau riemannien complet d'une 4-vari\'{e}t\'{e} compacte connexe est
conjugu\'{e} \`{a} un drapeau de Lie de flot homog\`{e}ne.
\end{corollary}

\begin{proof}
Soit $(M,\mathcal{F}_{1,}\mathcal{F}_{2},\mathcal{F}_{3})$ un tel drapeau.
Les deux derniers feuilletages \'{e}tant de Lie, il reste seulement \`{a}
montrer que le flot est de Lie.

- Si le drapeau est stable, le probl\`{e}me est r\'{e}solu.

- Si le flot est \`{a} feuilles ferm\'{e}es , le flot est un feuilletage de
Lie en cercles puisque la vari\'{e}t\'{e} basique du flot est le tore
hyperbolique $\mathbb{T}_{A}^{3}(trA\geq 2)$ (prop$.\ref{proposition2}$). Et
dans ce cas, $M=\mathbb{S}^{1}\times \mathbb{T}_{A}^{3}.$

- Si l'alg\`{e}bre de Lie structurale du flot est de dimension $un$, la vari%
\'{e}t\'{e} basique du flot \'{e}tant le tore $\mathbb{T}^{2},$ alors le
drapeau ($\mathcal{F}_{1}$, $\overline{\mathcal{F}_{1}}$) peut \^{e}tre compl%
\'{e}t\'{e} en un drapeau riemannien complet stable. Et ici aussi $\ M=%
\mathbb{T}^{2}\times _{A}\mathbb{T}^{2}$, \ o\`{u} $A$ appartient \`{a} la
famille des matrices d\'{e}crites plus haut.
\end{proof}

\subsection{Structure d'un drapeau riemannien complet}

Le second r\'{e}sultat fournit une obstruction \`{a} l'existence d'un
drapeau riemannien complet sur une \ vari\'{e}t\'{e} (non n\'{e}cessairement
compacte) \ et montre que les drapeaux riemannien complets ont des
propri\'{e}t\'{e}s assez riches et proches des drapeaux \ complets de
feuilletages de Lie homog\`{e}nes . Ce qui nous permet d'obtenir pour une
vari\'{e}t\'{e} compacte une caract\'{e}risation des drapeaux riemanniens
complets qui sont des drapeaux de feuilletages homog\`{e}nes.

\begin{theorem}
Soit $\mathcal{D=(F}_{1},...,\mathcal{F}_{n-1})$ un drapeau riemannien d'une 
$n-$ vari\'{e}t\'{e} connexe non n\'{e}cessairement compacte $M$. Alors on a
les propri\'{e}t\'{e}s suivantes.

1) La vari\'{e}t\'{e} $M$ est compl\`{e}tement parall\'{e}lisable et les
feuilletages du drapeau \ sont \`{a} la fois transversalement parall\'{e}%
lisable et transversalement int\'{e}grables. De fa\c{c}on pr\'{e}cise, il
existe un unique parall\'{e}lisme ($Y_{i})_{1\leq i\leq n}$ de $M$ dit 
\textbf{parall\'{e}lisme du drapeau} tel que:

- Pour tout $i\geq 1,Y_{i}$ est unitaire, tangent \`{a} $\mathcal{F}_{i},$et
oriente le flot qu'il d\'{e}finit

- Pour tous $1\leq i<j\leq n,$ $[Y_{i},Y_{j}]=k_{ij}Y_{i}$, o\`{u} les
fonctions $k_{ij}$, dites \textbf{fonctions de structure du drapeau}, sont ,
pour $i\geq 2,$ basiques pour le feuilletage $\mathcal{F}_{i-1}$.

2) Si la vari\'{e}t\'{e} est compacte, son r\'{e}v\^{e}tement universel est
isomorphe \`{a} $\mathbb{R}^{n}$et le drapeau relev\'{e} est un drapeau de
feuilletages simples transversalement parall\'{e}lisables. Et le drapeau$%
\mathcal{D}$ est un drapeau de feuilletages homog\`{e}nes si et seulement si
ses fonctions de structure sont constantes.
\end{theorem}

\begin{proof}
Pour un feuilletage riemannien ($M$, $\mathcal{F)}$ on d\'{e}signe par :

- $\mathcal{A}(M)$ l'anneau des fonctions num\'{e}riques $C^{\infty }$ sur $M
$

- $\Xi \mathcal{F}$ le$\mathcal{A}(M)-$ module des champs de vecteurs
tangents \`{a} $\mathcal{F}$,

- $\upsilon $($\mathcal{F}$) le fibr\'{e} normal de $\mathcal{F}$,

- $l(M$,$\mathcal{F}$) l'alg\`{e}bre de Lie des champs de vecteurs globaux
feuillet\'{e}s transverses,

- $(T\mathcal{F})^{\bot }$ le fibr\'{e} orthogonal de $T\mathcal{F}$
relativement \`{a} une m\'{e}trique quasi-fibr\'{e}e de $\mathcal{F}$ ,

- \ si $\ X\in \Xi \mathcal{F}$ , $\overline{<X>}$ \ le sous- fibr\'{e} de
rang\textit{\ un} qu'il engendre,

Soit $\mathcal{D=}$($\mathcal{F}_{1},...,\mathcal{F}_{p}$) un drapeau
riemannien sur la $n$-vari\'{e}t\'{e} riemannienne $(M$, $g),$ o\`{u} $g$
est une m\'{e}trique quasi-fibr\'{e}e du drapeau ( prop.\ref{prop.metrique
qf}) et $\mathcal{F}$ un feuilletage quelconque de $\mathcal{D}$. Soient $\ m
$=dim$\mathcal{F}$ ( $i.e.$ $\mathcal{F}$ =$\mathcal{F}_{m}$ ) . Alors pour
tout $i,m\leq i\leq n-2$ , les sous-fibr\'{e}s ($T\mathcal{F}_{i})^{\bot
}\cap T\mathcal{F}_{i+1}$ et ($T\mathcal{F}$ $_{n-1}$)$^{\bot }$ sont de
rang $un$ ( l'orthogonalit\'{e} \'{e}tant bien s\^{u}r entendue dans le sens
de la m\'{e}trique quasi-fibr\'{e}e du drapeau). Soient respectivement $%
Y_{i+1}$ et $Y_{n}$ les champs de vecteurs unitaires qui les engendrent et
qui orientent chacun des flots associ\'{e}s.

D'apr\`{e}s ce qui pr\'{e}c\`{e}de, puisque $g$ est quasi-fibr\'{e}e pour $%
\mathcal{F}$, il est clair que:

1) $Y_{m+1}$,..,$Y_{n}$ sont $(n-m)$ champs globaux lin\'{e}airement ind\'{e}%
pendants en chaque point de $M$, formant une base de $\upsilon $($\mathcal{F}%
_{m}$)$\cong $ (T$\mathcal{F}_{m}$)$^{\bot }$ .

2) Pour tout $k,0\leq k\leq n-m-1$%
\begin{equation*}
T\mathcal{F}_{m+k+1}=T\mathcal{F}_{m+k}\oplus \overline{<Y_{m+k+1}>}=T%
\mathcal{F}_{m}\oplus (\underset{m+1\leq j\leq m+k+1}{\oplus \overline{%
<Y_{j}>})}
\end{equation*}

3) Pour tout $i\geq m+1$ et pour tout $X\in \Xi \mathcal{F}$, $Y_{i}$ \'{e}%
tant unitaire alors on a: 
\begin{equation*}
0=Xg(Y_{i},Y_{i})=2g([X,Y_{i}],Y_{i}).
\end{equation*}

Il r\'{e}sulte de 2) et 3) que pour tout $i\geq m+1$, [ $X$,$Y_{i}$] $\in
\Xi \mathcal{F}_{i}$ et est orthogonal \`{a} $Y_{i}$ ; ce qui implique que $%
[X,Y_{i}]$ est tangent \`{a} $\mathcal{F}_{i-1}$ puisque \newline
$T\mathcal{F}_{i}=T\mathcal{F}_{i-1}\oplus \overline{<Y_{i}>}.$ Ensuite pour
tous $\ i,$ $j,$ $\ m+1\leq i<j\leq n$, on a 
\begin{equation*}
0=Xg(Y_{i},Y_{j})=g([X,Y_{i}],Y_{j})+g(Y_{i},[X,Y_{j}])=g([X,Y_{i}],Y_{i})
\end{equation*}
car [ $X$,$Y_{i}$]$\in \Xi \mathcal{F}_{i-1}$ est orthogonal \`{a} $Y_{j}$
puisque $Y_{j}$ $\bot T\mathcal{F}_{j-1}$ et $T\mathcal{F}$ $_{i-1}$ $%
\subset T\mathcal{F}_{j-1}.$

Ceci montre que $[X,Y_{i}]$ est tangent \`{a} $\mathcal{F}_{i-1}$ et est
orthogonal \`{a} $Y_{m+1}$,..,$Y_{n}$. Comme $T\mathcal{F}$ $_{i-1}$= $T%
\mathcal{F}_{m}\oplus (\underset{m+1\leq j\leq i-1}{\oplus }\overline{<Y_{k}>%
})$ \ alors $\ [X,Y_{i}]$ est tangent \`{a} $\mathcal{F}_{m}$, $i.e$. pour
tout $i\geq m+1$,$\ Y_{j}\in l$($M$,$\mathcal{F}$). Ce qui assure que $%
\mathcal{F}$ est un feuilletage transversalement parall\'{e}lisable; par
suite le drapeau est parall\'{e}lisable.

Pour $m=1,$ soit $Y_{1}$ \ un champ unitaire complet tangent au flot $%
\mathcal{F}_{1},(Y_{2},..,Y_{n})$ le parall\'{e}lisme tranverse de $\mathcal{%
F}_{1}$ , alors par construction, $(Y_{1},..,Y_{n})$ est un parall\'{e}lisme
de $TM$ \ =$T\mathcal{F}_{1}\oplus (\underset{2\leq j\leq n}{\oplus 
\overline{<Y_{j}>})}=(\underset{1\leq j\leq n}{\oplus \overline{<Y_{j}>})}$
et $(Y_{i+1},..,Y_{n})$ un parall\'{e}lisme de $\mathcal{F}_{i}.$ Ce parall%
\'{e}lisme est alors tel que :\newline
pour tous $i,j,1\leq i\leq j\leq n,[Y_{i},Y_{j}]\in \Xi \mathcal{F}_{i}$ et
pour tout $k<i,$ 
\begin{equation*}
g([Y_{i},Y_{j}],Y_{k})=Y_{j}g(Y_{k},Y_{i})-g([Y_{k},Y_{j}],Y_{i})=0
\end{equation*}
Il en r\'{e}sulte que $[Y_{i},Y_{j}]\in \overline{<Y_{1},..,Y_{i-1}>}^{\bot
}\cap \Xi \mathcal{F}_{i}=\overline{<Y_{i}>}$ , donc pour tous $i,j,$ $1\leq
i<j\leq n$ on a 
\begin{equation*}
\lbrack Y_{i},Y_{j}]=k_{ij}Y_{i}\text{, (**)}
\end{equation*}
et alors la vari\'{e}t\'{e} $M$ est compl\`{e}tement parall\'{e}lisable.

Pour $l<i<j,$ la relation de Jacobi appliqu\'{e}e au triplet ( $Y_{l}$,$Y_{i}
$, $Y_{j})$ et l'identit\'{e}(**) permettent de voir que les fonctions de
structure $k_{ij}$ sont basiques pour le feuilletage $\mathcal{F}_{i-1}.$%
L'identit\'{e} (**) montre \'{e}galement que le $\mathcal{A}(M)-\Xi $ $%
\upsilon $($\mathcal{F}_{i}$) $\ $des champs normaux \`{a} $\mathcal{F}_{i}$
est un module libre de base $Y_{i+1},..,Y_{n},$ et que les distributions $%
\upsilon $($\mathcal{F}_{i}$)$\approxeq (T\mathcal{F}_{i})^{\bot }$ sont int%
\'{e}grables. On notera pour la suite, $\mathcal{F}_{i}{}^{\bot }$ le
feuilletage transverse ainsi d\'{e}fini.

2) \ Supposons maintenant que la vari\'{e}t\'{e} de base $M$ est compacte.

a) Dans ce cas la m\'{e}trique $g$ et sa relev\'{e}e$\widetilde{g}$ sur le
rev\^{e}tement universel de $M$ \'{e}tant compl\`{e}tes, et $\mathcal{F}_{i}$
transversalement int\'{e}grable, le th\'{e}or\`{e}me de d\'{e}composition \
de Blumenthal- Hebda \cite{BH}assure que le rev\^{e}tement universel $%
\widetilde{M\text{ }}$de $M$ s'identifie au produit $\widetilde{F_{i}}\times 
\widetilde{F_{i}}^{\bot },$ o\`{u} $\widetilde{F_{i}}$ est une feuille du
feuilletage relev\'{e}$\widetilde{\mathcal{F}_{i}}$ de $\mathcal{F}_{i}$ et $%
\widetilde{F_{i}}^{\bot }$une feuille du feuilletage relev\'{e} $\widetilde{%
\mathcal{F}_{i}^{\bot }}$de $\mathcal{F}_{i}^{\bot },$et que la projection
sur le second facteur est une submersion riemannienne qui d\'{e}finit le
feuilletage relev\'{e} $\widetilde{\mathcal{F}_{i}}$.Il en r\'{e}sulte comme 
$\widetilde{F_{1}}$ est un rev\^{e}tement simplement connexe de la feuille $%
F_{1}$ de $\mathcal{F}_{1}$ projet\'{e}e de $\widetilde{F_{1}}$ , que $%
\widetilde{M\text{ }}\approxeq \mathbb{R}\times \widetilde{F_{1}}^{\bot }.$
Comme le drapeau relev\'{e} sur ($\widetilde{M\text{ }},\widetilde{g}$ $)$
est aussi un drapeau parall\'{e}lisable et comme la projection sur le second
facteur est la fibration basique du flot relev\'{e}, alors ledrapeau relev%
\'{e} se projette sur $\widetilde{F_{1}}^{\bot }$en un drapeau riemannien
complet d'une vari\'{e}t\'{e} compl\`{e}te connexe et simplement connexe et
le m\^{e}me th\'{e}or\`{e}me de d\'{e}composition \ de Blumenthal- Hebda
assure que $\widetilde{F_{1}}^{\bot }\approxeq \mathbb{R}\times \widetilde{%
F_{2}}^{\bot }$ et on a $\widetilde{M\text{ }}\approxeq \mathbb{R}^{2}\times 
\widetilde{F_{2}}^{\bot }.$ De proche en proche on obtient que $\widetilde{M%
\text{ }}\approxeq \mathbb{R}^{n-1}\times \widetilde{F_{n-1}}^{\bot
}\approxeq \mathbb{R}^{n}.$ On notera de m\^{e}me que le rev\^{e}tement
universel des feuilles du feuilletage $\mathcal{F}_{i}$ du drapeau sont
toutes isomorphes \`{a} $\mathbb{R}^{i}.$

b) - Si $(M,\mathcal{D)}$ est un drapeau de feuilletages homog\`{e}nes, le th%
\'{e}or\`{e}me de rel\`{e}vement d'un drapeau de Lie\cite{DIA1}, rel\`{e}ve $%
\mathcal{D}$ en un drapeau de feuilletages homog\`{e}nes simples sur un
groupe de Lie $G$ connexe et simplement connexe rev\^{e}tement universel de $%
M.$ Comme l'action \`{a} gauche d'\'{e}l\'{e}ments de $G$ est une isom\'{e}%
trie et que les feuilletages rel\'{e}v\'{e}s sont homog\`{e}nes, il est
clair que le parall\'{e}lisme de ce drapeau rel\'{e}v\'{e} est compos\'{e}
de champs invariants \`{a} gauche.Si les $C_{ij}^{k}$ sont les constantes de
structure de l'alg\`{e}bre de Lie de $G$ dans la base d\'{e}finie par ce
parall\'{e}lisme, il est alors clair que les fonctions de structure de $%
\widetilde{k_{ij}}$ de ce drapeau relev\'{e} sont alors constantes, soit $%
\widetilde{k_{ij}}=C_{ij}^{j}$ . Comme la projection sur $M$ de ce parall%
\'{e}lisme est le parall\'{e}lisme du drapeau $(M,\mathcal{D)}$ avec les m%
\^{e}mes fonctions de structure, alors on a le r\'{e}sultat.

- R\'{e}ciproquement \ si les fonctions de structure $k_{ij}$ du parall\'{e}%
lisme ($Y_{i})_{1\leq i\leq n}$ d'un drapeau riemannien complet $(M,\mathcal{%
D)}$ sont constantes, alors le parall\'{e}lisme rel\'{e}v\'{e} ($\widetilde{Y%
}_{i})_{1\leq i\leq n}$ sur le rev\^{e}tement universel $\widetilde{M\text{ }%
}$de $M$, d\'{e}finit une alg\`{e}bre de Lie r\'{e}elle de dimension $n,$ de
groupe de Lie connexe et simplement connexe $G.$ On peut supposer , \`{a} un
diff\'{e}omorphisme pr\`{e}s, que $\widetilde{M}$ =$G$ et que $M$ est un
espace homog\`{e}ne. Puisque le feuilletage rel\'{e}v\'{e} de tout
feuilletage $\mathcal{F}_{i}$ du drapeau est, comme on le voit , d\'{e}fini
par les orbites du sous-groupe distingu\'{e} correspondant \`{a} l'id\'{e}al
engendr\'{e} par $\widetilde{Y_{1}},..,$ $\widetilde{Y_{i}}$ (**), alors le
feuilletage $\mathcal{F}_{i}$ est un feuilletage homog\`{e}ne.
\end{proof}

\bibliographystyle{AABBRV}
\bibliography{,,,,,,,,,,,,acompat}

\end{document}